\documentclass[twocolumn]{autart}    
\usepackage{savesym}
\savesymbol{AND}
\usepackage{graphicx}
\usepackage{amsmath}
\usepackage{amssymb}
\usepackage{caption}
\usepackage{subcaption}
\usepackage{epstopdf}
\usepackage{algorithmic}
\usepackage{algorithm}
\usepackage{bm}
\usepackage{cite}

\DeclareMathOperator*{\argmin}{arg\,min}

\DeclareMathOperator{\proj}{proj}
\DeclareMathOperator{\E}{\mathbb{E}}

\DeclareMathOperator{\rounding}{rounding}
\DeclareMathOperator{\card}{card}

\DeclareMathOperator{\conv}{conv}

\newcommand{\RR}{\mathbb{R}}
\renewcommand{\SS}{\mathcal{S}}
\newcommand{\Rt}{\mathcal{R}}
\newcommand{\xx}{\mathbf{x}}
\newcommand{\yy}{\mathbf{y}}
\renewcommand{\gg}{\mathbf{g}}
\newcommand{\hf}{\hat{f}}
\newcommand{\ft}{\tilde{f}}
\newcommand{\ftg}{\ft^\text{g}}

\newcommand{\tSs}{\tilde{S}^\star}
\newcommand{\cc}{\bm{\chi}}
\newcommand{\R}{\texttt{R}}

\begin{document}

\begin{frontmatter}

\title{Online Dynamic Submodular Optimization\thanksref{footnoteinfo}} 
\thanks[footnoteinfo]{This paper was not presented at any IFAC 
meeting. Corresponding author A.~Lesage-Landry. Tel. +1-514-340-4711, Ext. 2442. 
Fax +1-514-340-4078.}

\author[Poly]{Antoine Lesage-Landry}\ead{antoine.lesage-landry@polymtl.ca},    
\author[Poly]{Julien Pallage}\ead{julien.pallage@polymtl.ca}              

\address[Poly]{Department of Electrical Engineering, Polytechnique Montreal, GERAD \& Mila, Montreal, QC, Canada, H3T 1J4.}  
          
\begin{abstract}                          
We propose new algorithms with provable performance for online binary optimization subject to general constraints and in dynamic settings. We consider the subset of problems in which the objective function is submodular. We propose the online submodular greedy algorithm (\texttt{OSGA}) which solves to optimality an approximation of the previous round loss function to avoid the NP-hardness of the original problem. We extend \texttt{OSGA} to a generic approximation function. We show that \texttt{OSGA} has a dynamic regret bound similar to the tightest bounds in online convex optimization with respect to the time horizon and the cumulative round optimum variation. For instances where no approximation exists or a computationally simpler implementation is desired, we design the online submodular projected gradient descent (\texttt{OSPGD}) by leveraging the Lovász extension. We obtain a regret bound that is akin to the conventional online gradient descent (\texttt{OGD}). Finally, we numerically test our algorithms in two power system applications: fast-timescale demand response and real-time distribution network reconfiguration.
\end{abstract}

\end{frontmatter}

\section{Introduction}
Online dynamic decision-making aims to consecutively provide decisions to minimize each round's objective function while relying only on the outcome of previous rounds.
The objective function is considered to be time-varying and decisions are made at each discretized time instance. Moreover, the objective function is assumed to be unknown at the time decisions have to be made. The online convex optimization (OCO)~\cite{zinkevich2003online,hazan2012online,shalev2012online} framework assumes a convex objective function and a convex and compact decision set. Provable performance guarantees can be established under some additional assumptions, e.g., boundedness of the objective function and its gradient~\cite{zinkevich2003online}, or the cumulative difference in round optima computed in hindsight~\cite{zinkevich2003online,mokhtari2016online}.

Online optimization is an appealing framework for real-time decision-making problems because it uses computationally efficient and scalable updates and provides performance guarantees. For example, it is used in the context of moving target tracking~\cite{lesage2020second,pun2020dynamic}, resource allocation in data centers~\cite{chen2017online,cao2018virtual}, portfolio selection~\cite{hazan2016introduction}, internet-of-things~\cite{chen2018bandit}, dynamic pricing in power systems~\cite{kim2016online}, or renewable intermittency mitigation~\cite{lesage2018setpoint,lesage2021online}.

Online binary optimization~\cite{jegelka2011online,hazan2012online,lesage2021online} considers a subset of problems in which the feasible set is the intersection of an application-specific constraint set and the binary set $\left\{0,1 \right\}^n$, where $n \in \mathbb{N}$ is the decision variable's dimension. 
{ It is motivated by real-time network flow, routing, scheduling, and knapsack problems which appear in engineering fields like telecommunications, logistics and operations, and electric power systems.
To tackle efficiently constrained, non-linear online binary optimization problems, we further assume that the objective function is submodular. Specifically, we consider online dynamic submodular optimization for which the objective is to provide the round optimal binary decisions. 
We propose two types of algorithms: (i) greedy approaches that solve approximations of the previous round's objective function and (ii) a projected gradient-based approach using the continuous and convex Lov\'asz extension of submodular functions. For all algorithms, we provide a performance analysis based on the dynamic regret. 
The regret bounds are shown to be sublinear in the number of rounds under different conditions on the variation between round optima computed in hindsight. Under these assumptions, the time-averaged dynamic regret vanishes as the time horizon increases and is, therefore, Hannan-consistent~\cite{hazan2016introduction}.

\paragraph*{Related work}
We now review the relevant literature on online non-linear binary optimization. Linearity simplifies considerably the problem as argued by~\cite{jegelka2011online} and, for this reason, is not considered. Examples of online binary optimization approaches for linear problems include~\cite{kalai2005efficient,koolen2010hedging}. In~\cite{hazan2012online}, the authors first studied the online submodular optimization problem. They only considered the static setting in which the decisions are benchmarked with a single static decision computed in hindsight. This is referred to as a static regret analysis~\cite{zinkevich2003online}. They further restricted their analysis to unconstrained problems. Reference~\cite{jegelka2011online} then proposed approaches to integrate constraints within online submodular optimization. They also limited their analysis to the static setting. In both cases, greedy and projected gradient-based approaches are proposed. In this work, we present a dynamic regret analysis for all our approaches which in turns provides a performance guarantee with respect to the round optimum. This latter aspect is important in an engineering setting because one wants to achieve optimality at each round, e.g., to track a time-varying setpoint. 
In the dynamic setting,~\cite{lesage2021online} used randomization and online convex optimization to solve problems with convex objective functions, i.e., convex with respect to the convex hull of the decision set. However,~\cite{lesage2021online} do not admit other than binary constraints and the dynamic regret analysis does not hold asymptotically. 

In the power system literature, several approaches based on time-varying optimization have been proposed to deal with binary decision variables. References~\cite{bernstein2016real,bernstein2019real} applied the error diffusion algorithm to obtain binary decisions from continuous decisions computed via the relaxed problem. In~\cite{zhou2019online}, randomization is used to convert continuous decisions to binary ones. This body of literature does not compare the round minima with the algorithm's decisions like online optimization does using the dynamic regret. Specifically, we make the following contributions:
\begin{itemize}
  \item We propose two algorithms for online dynamic submodular optimization. Under the submodularity assumption, we provide online constrained binary optimization algorithms with provable performance guarantees in dynamic settings which hold for the first time when subject to constraints and/or any time horizon. We extend the static regret analysis of~\cite{jegelka2011online,hazan2012online} and establish conditions under which our algorithms lead to a sublinear dynamic regret bound in the number of rounds and (tractable) round optimum variation.
  \item We formulate a greedy algorithm that solves a $\beta$-approximation of the previous round's objective function. When this approximation is not available, we show that a generic approximation can be used with limited impact on the performance bound.
  \item We provide a computationally efficient and scalable algorithm for fast timescale online optimization problems which only performs a single project gradient descent step on the Lov\'asz extension of the objective function.
  \item We numerically evaluate the performance of our approaches in power system examples. First, we use the projected gradient-descent update to dispatch demand response resources for frequency regulation. Second, we apply the greedy update to real-time network configuration where line switches can be controlled (\textsc{on}/\textsc{off}) to minimize the active power losses while spanning a radial network.
\end{itemize}

Next, we provide background on online optimization and submodularity and introduce our notation. Greedy and projected gradient descent-based approaches are analyzed in Sections~\ref{sec:greedy}~and~\ref{sec:projected}, respectively. Numerical examples showcasing our approaches in power system applications are presented in Section~\ref{sec:appli}. Conclusions and future work are provided in Section~\ref{sec:conclu}.

\section{Preliminaries}
In this section, we introduce our notation and the online optimization setting, and provide the relevant background on submodular functions.

\subsection{Online optimization}
In online optimization, a round-dependent objective function must be minimized at each round $t \in \left\{1,2, \ldots, T \right\}$, where $T \in \mathbb{N}$ is the time horizon. In this setting, the objective function is assumed to be observed only after the decision maker has implemented the round's decision, which must be provided on a fast timescale.

We consider a subset of online binary optimization problems with the base set $V = \left\{1,2, 3, \ldots, n \right\}$, $n \in \mathbb{N}$, in which the objective function is assumed to be submodular. Let the power set $2^V$ represent the set of all possible decisions. In each round $t \in \left\{1,2, \ldots, T \right\}$, a decision $S \in 2^V$ must be made. The problem takes the form:
\begin{equation}
\min_{S \in \SS} f_t(S), \label{eq:prob_submodular}
\end{equation}
where $f_t: 2^V \mapsto \mathbb{R}$ is a submodular set function, $\SS \subseteq 2^V$ is the feasible set, i.e., the set that expresses the problem's constraints, and $t \in \left\{1,2, \ldots, T \right\}$. 

As noted by~\cite{jegelka2011online}, at time $t$,~\eqref{eq:prob_submodular} is NP-hard if $\SS \neq 2^V$. Because no offline optimization algorithm can solve~\eqref{eq:prob_submodular} given $f_t$ in polynomial time, we benchmark the decisions provided by our online optimization algorithm with an offline $\alpha$-approximation algorithm~\cite{jegelka2011online}. Let $S^\star_t \in \argmin_{S \subseteq \SS} f_t(S)$. 
An $\alpha$-approximation algorithm provides a solution $S_t^*$ such that $f_t(S_t^*) \leq \alpha f_t(S_t^\star)$. For common submodular minimization problems like minimum spanning tree~\cite{goel2009approximability} or edge cover~\cite{iwata2009submodular}, $\alpha$ values are related to their graph structure~\cite{jegelka2011online}.
Building on~\cite{jegelka2011online}, we define the \emph{dynamic} $\alpha$-regret to characterize the performance of our online optimization approaches.

\begin{defn}
The dynamic $\alpha$-regret $\textup{\R}^\textup{\text{d}}_\alpha(T)$ over a time horizon $T$ is:
\[
\textup{\R}^{\textup{d}}_\alpha(T) = \sum_{t=1}^T \left( f_t(S_t) - \alpha f_t(S_t^\star) \right),
\]
where $S_t$ is the decision provided by the online optimization algorithm at round $t$.
\end{defn}
We note that the special case $\SS = 2^V$ can be solved to optimality in polynomial time. At this time, we let $\alpha = 1$ and we retrieve the standard dynamic regret definition from OCO~\cite{zinkevich2003online}. This fact is used to specialize our results.

The dynamic $\alpha$-regret defers from the $\alpha$-regret employed in~\cite{jegelka2011online} because it uses as comparators (second term of the sum) to the algorithm's decision (first term), the round optima instead of the best fixed decision in hindsight. The power system applications, later discussed in Section~\ref{sec:appli}, motivate the use of a framework that targets round optimal decisions instead of an averaged, static decision. Dynamic regret bounds are given in terms of the cumulative round optimum variation $V_T$ or a derivative of it~\cite{zinkevich2003online,hall2015online}. This term is used as a complexity measure in dynamic problem~\cite{mokhtari2016online}.

We conclude by adapting the definition of $V_T$ to online set function optimization from which online submodular optimization is a subset of. Let $\bm{\chi}_A \in \{0,1\}^n$ where $\bm{\chi}_A$'s $i^\text{th}$ component is one if and only if $i \in A$ and zero otherwise be the characteristic vector of the set $A \in 2^V$. 
Consider the online binary optimization problem counterpart of~\eqref{eq:prob_submodular}:
\begin{equation*}
\min_{\yy \in {\mathcal{Y} \cap \left\{0,1\right\}^n}} f_t^\text{b}(\yy), \label{eq:prob_combinatorial}
\end{equation*}
where $\mathcal{Y} \subseteq \RR^n$ is the constraint set and $f^\text{b}: \left\{0,1\right\}^n \mapsto \RR$ is the objective function. Let $\yy_{t}^\star \in \argmin_{\yy \in \mathcal{Y} \cap \left\{0,1\right\}^n} f^\text{b}\left(\yy\right)$. The cumulative variation term $V_T$, as in standard online (convex) optimization, is~\cite{zinkevich2003online,lesage2021online}:
\begin{align*}
V_T &= \sum^T_{t=2} \left\| {\yy}_{t}^\star - \yy_{t-1}^\star \right\|_2 = \sum^T_{t=2} \left\| \cc_{S_{t}^\star} - \cc_{S_{t-1}^\star} \right\|_2,
\end{align*}
where $S_{t}^\star \in 2^V$ is the subset of components of $\yy_{t}^\star\in\left\{0,1\right\}^n$ with value one. Adapting $V_T$ to set-valued objective functions, we, therefore, obtain:
\begin{equation*}
V_T = \sum^T_{t=2} \sqrt{\card\left( S_{t}^\star \ominus S_{t-1}^\star\right)}, \label{eq:V_T}
\end{equation*}
where $\ominus$ is the symmetric difference or disjunctive union of two sets. When the context requires it, we will introduce alternative $V_T$ definitions, e.g., when the optima are defined from function approximations.

For future comparison, previous work in online optimization in static settings~\cite{jegelka2011online,hazan2012online} establish $O( \sqrt{T})$ regret and $\alpha$-regret bounds, respectively, where $S^\star_t = S^\star \in \argmin_{S\in\SS} \sum_{t=1}^T f_t (S) \ \forall t$ is used as the static comparator. In the dynamic, unconstrained setting,~\cite{lesage2021online} obtained an $O\left(T^{\epsilon} + V_T \right)$, where $\epsilon \in (0,1)$ regret bound for a finite time horizon $T$.

\subsection{Submodularity}
\label{ssec:sub}
A function $f_t:2^V \mapsto \RR$ is submodular if it exhibits the diminishing marginal return property~\cite{jegelka2011online}, i.e., if
$f_t(A \cup \left\{ i \right\}) - f_t(A) \geq f_t(B \cup \left\{ i \right\}) - f_t(B)$,
for all $A \subseteq B \subseteq V$ and $i \in V$. The Lov\'asz extension $\hf: [0,1]^n \mapsto \RR$ of a function $f_t$ can be defined as:
\[
\hf_t(\xx) = \sum_{i=1}^n x_i \left( f_t(\left\{1, 2, \ldots, i \right\}) - f_t(\left\{1, 2, \ldots, i-1 \right\}) \right),
\]
where $x_i$ is the $i^\text{th}$ largest component of $\xx \in [0,1]^n$, $\{0\}\equiv\emptyset$, and $f_t(\emptyset) \equiv 0$~\cite{bach2013learning}. Lastly, we will make use of two important properties of the Lov\'asz extension: (i) $\hf_t$ is convex if and only if $f_t$ is submodular and (ii) $\hf_t\left( \bm{\chi}_A \right) = f_t(A)$ for submodular functions.

A subgradient of $\hf_t$ at a point $\xx \in \conv \SS$ can be computed using only evaluations of the original, submodular function $f_t$. Let $\pi: \left[0, 1 \right]^n \times V \mapsto V$ be a function where $\pi(\xx, i) = j$ is such that the $i^\text{th}$ largest component of $\xx$ is $x_j$. Let $\partial \hf_t(\xx)$ be the subgradient set of $\hf$ at $\xx$. Then, we have the following definition for $\gg_t \in \partial \hf_t(\xx)$:
\begin{equation}
\begin{aligned}
\gg_t = \sum_{i=1}^n &\left( f_t(\left\{1, 2, \ldots, i \right\})\right. \\
&\quad\left.- f_t(\left\{1,  2,\ldots, i-1 \right\}) \right) \bm{\chi}_{\{ \pi(\xx, i) \}} .
\end{aligned}
\label{eq:gradient}
\end{equation}

Finally, a rounding algorithm, $\rounding_\SS:\left[0,1\right]^n \mapsto \SS $, can be employed to convert the Lova\'sz extension's continuous input $\xx \in \left[0,1\right]^n$ to the corresponding set of the original, set function $f_t$. 
For example, in Section~\ref{sec:appli} we will use $\rounding_{2^V}(\xx): \left[0,1\right]^n \mapsto 2^V$, a standard rounding map for unconstrained problems defined as follows. Let $\xx \in \left[0,1\right]^n$, then $\rounding_{2^V}(\xx) = S$ where $S = \left\{\left. i \in V \right| x_i \geq p \right\}$ with $p \sim \text{Uniform}[0,1]$. Note that we get $\hf_t(\xx) = \E\left[ f_t(S)\right]$~\cite{hazan2012online}. Alternatively, for some types of feasible sets $\SS$~\cite{iwata2009submodular,jegelka2011approximation}, rounding algorithms can be characterized by their approximation guarantee~\cite{jegelka2011online}. For example, a rounding technique $\rounding_\SS$ with approximation guarantee $\alpha$ is such that $\alpha \hf_t(\xx_t) \geq f_t(S_t) = \hf_t(\cc_{S_t})$ for $S_t = \rounding_\SS(\xx_t)$, where $S_t \in \SS$ and $\xx_t \in \conv \SS$.

\section{Greedy approaches}
\label{sec:greedy}
We now propose greedy approaches for online binary optimization. These approaches are based on the previous round's objective function and an approximation that renders the submodular problem tractable. We first consider the following function approximation.

\begin{defn}[$\beta$-approximation function~\cite{jegelka2011online}]
The function $\tilde{f}_t: 2^V \mapsto \RR$ is a $\beta$-approximation of $f_t$ if it satisfies the following conditions:
\begin{enumerate}
  \item $f_t(S) \leq \ft_{t}(S) \leq \beta f_t(S)$ for $\beta \geq 1$ and all $S \subseteq V$;
  \item $\min_{S \subseteq \SS} \ft_{t}(S)$ can be solved to optimality in polynomial time. \label{def:beta_approx_poly}
\end{enumerate}
\label{def:beta_approx}
\end{defn}

Examples of approximations compatible with Definition~\ref{def:beta_approx} are provided in~\cite[Section 2]{jegelka2011online}. Note that we require the milder condition where $\min_{S \subseteq \SS} \ft_t(S)$ needs to be tractable contrarily to~\cite{jegelka2011online} which impose the condition on $\min_{S \subseteq \SS} \sum_{i=1}^t \tilde{f}_i(S)$ at round $t$.
Based on Definition~\ref{def:beta_approx},~\eqref{eq:prob_submodular} can be tackled in an online fashion using the following update:
\begin{equation}
S_t \in \argmin_{S \in \SS} \tilde{f}_{t-1}(S), \label{eq:greedy_update}
\end{equation}
where $\ft$ is a $\beta$-approximation function of $f$. We refer to an algorithm implementing~\eqref{eq:greedy_update} as online submodular greedy algorithm (\texttt{OSGA}). By Definition~\ref{def:beta_approx}$-$\eqref{def:beta_approx_poly}, \eqref{eq:greedy_update} can be efficiently solved to optimality.

For the next results, we make the following assumptions. 
\begin{assum}
Let $f_t$ be a bounded function over the set $\mathcal{S}$, i.e., there exists $M \in \RR_{>0}$ such that $|f_t(S)|\leq M$ for all $t=1,2,\ldots, T$ and $S \in \SS$.
\label{ass:bounded}
\end{assum}
\begin{assum}
The set value function $f_t$ is such that
$\left| f_t(S_1) - f_t(S_2) \right| \leq L \card\left(S_1 \ominus S_2 \right),$
for all $S_1, S_2 \subseteq \mathcal{V}$ and $0 < L < +\infty$.
\label{ass:pseudo_L}
\end{assum}
In other words, we assume a Lipschitz continuity-like property for set functions. Assumption~\ref{ass:pseudo_L} holds for any submodular function if Assumption~\ref{ass:bounded} does, e.g., the generic approximation defined below~\cite{goemans2009approximating} or the $\beta$-approximation function for minimum spanning tree with submodular loss function $h(S)$, $\tilde{h}(S) = \sum_{i \in S} h(i)$~\cite{goel2009approximability}.

For the regret analysis, we let $\tilde{S}^\star_t \in \argmin_{S \in \SS} \tilde{f}_t(S)$ where $\tilde{f}_t$ is a $\beta$-approximation of $f_t$. We redefine the cumulative variation of the optima as $\tilde{V}_T = \sum_{t=2}^T\left\| \chi_{\tSs_{t}} - \chi_{\tSs_{t-1}}\right\|_2$. This definition is similar to the one used in standard dynamic online convex optimization~\cite{zinkevich2003online,hall2015online} and has the advantage of being a function of efficiently obtainable optima. We recall that our objective is to establish a dynamic $\alpha$-regret bound for an online optimization algorithm tackling~\eqref{eq:prob_submodular}, i.e., to achieve performance similar to an offline $\alpha$-approximation algorithm used consecutively. Providing a dynamic $\alpha$-regret bound in term of $\tilde{V}_T$ is in submodular line with this objective because it can be effectively characterized as opposed to $V_T$ which requires solving a sequence of $T$ NP-hard problems. The regret analysis of update~\eqref{eq:greedy_update} is provided in Theorem~\ref{eq:thm_greedy}.

\begin{thm}
Suppose $\ft_t$ is a $\beta$-approximation of $f_t$ such that $\ft_t$ satisfies Assumptions~\ref{ass:bounded} and~\ref{ass:pseudo_L}. 
If $\alpha \geq \beta$, then the $\alpha$-regret of \texttt{OSGA} is bounded by:
\[
\textup{\R}^{\textup{d}}_\alpha(T) \leq \frac{\alpha L}{\beta} \sum^T_{t=2} \sqrt{\card\left( \tSs_{t} \ominus \tSs_{t-1}\right)} =  \frac{\alpha L}{\beta} \tilde{V}_T.
\]
If $\tilde{V}_T$ is sublinear, then so is the $\alpha$-regret.
\label{eq:thm_greedy}
\end{thm}
\begin{pf}
We bound the $\alpha$-regret using Definition~\ref{def:beta_approx} to obtain
\begin{equation}
\R^\text{d}_\alpha(T) \leq \sum_{t=1}^T \tilde{f}_t(S_t) - \frac{\alpha}{\beta} \tilde{f}_t(S^\star_t).\label{eq:inter_tilde}
\end{equation}
We observe that $\tilde{f}_t(S_t) = \tilde{f}_t(\tilde{S}_{t-1}^\star)$ because of~\eqref{eq:greedy_update} and $\tilde{S}_t^\star \in \argmin_{S \in \SS} \tilde{f}_t(S)$. Thus, we can rewrite~\eqref{eq:inter_tilde} as
\begin{align}
\R^\text{d}_\alpha(T) & \leq \sum_{t=1}^T \tilde{f}_t(\tilde{S}_{t-1}^\star) - \frac{\alpha}{\beta} \tilde{f}_t(\tilde{S}^\star_t), \label{eq:inter_alphabeta}
\end{align}
where we also used the definition of $\tilde{S}^\star_t$. 
By assumption, $\alpha \geq \beta$ and factoring out $\alpha/\beta$ of~\eqref{eq:inter_alphabeta}'s sum yields an upper bound. Then, by Assumption~\ref{ass:pseudo_L}, we obtain
\begin{align*}
\R^\text{d}_\alpha(T) & \leq \frac{\alpha L}{\beta} \sum_{t=1}^T \sqrt{\card \left( \tilde{S}_t^\star \ominus \tilde{S}_{t-1}^\star \right)},
\end{align*}
and we have completed the proof. $\hfill \qed$
\end{pf} 
We remark that contrarily to~\cite[Theorem 2]{jegelka2011online}, the approximation factor $\alpha$ does not need to be known to run the algorithm. Given a $\beta$-approximation of $f_t$, \texttt{OSGA} leads to an $O(\tilde{V}_T)$ regret bound that resembles the tightest dynamic bound in standard OCO~\cite{mokhtari2016online}, i.e., $O(1+V_T)$, with respect to (w.r.t.)~$T$ and the variation of tractable optima. Note that: (i)} this latter work requires strong convexity and (ii) our results uses the $\alpha$-approximation algorithm's solutions as comparators in the regret. Theorem~\ref{eq:thm_greedy}'s bound also improves on~\cite{lesage2021online}'s expected bound because it is only a function of the cumulative variation of obtainable optima and holds asymptotically. 

Recall that unconstrained submodular minimization problem can be solved to optimality efficiently. Hence,
for the special case where $\mathcal{S} = 2^V$, i.e., when~\eqref{eq:prob_submodular} is an unconstrained submodular problem, \texttt{OSGA} can be directly applied to the previous round loss function. This application leads to the following regret bound.

\begin{cor}
If $\SS = 2^V$ and $\alpha = \beta = 1$, then \texttt{OSGA}'s update reduces to
\begin{equation}
S_t \in \argmin_{S \in 2^V} f_{t-1}(S), \label{eq:greedy_update_V}
\end{equation}
and leads to:
\[
\textup{\R}^\textup{\text{d}(T)} \leq L \sum_{t=2}^T \sqrt{\card \left( {S_{t}^\star} \ominus {S_{t-1}^\star}\right)} < O\left( V_T \right).
\]
\end{cor}
\begin{pf}
The proof follows from Theorem~\ref{eq:thm_greedy} where (i) the regret is considered instead of the $\alpha$-regret and (ii) $f_t$ is used directly instead of $\tilde{f}_t$, the $\beta$-approximation because~\eqref{eq:greedy_update_V} can be solved efficiently. $\hfill \qed$
\end{pf}

For some problem instances, finding an approximation that satisfies both Definition~\ref{def:beta_approx} and Assumption~\ref{ass:pseudo_L} is difficult. Alternatively, the generic approximation for submodular functions provided in Definition~\ref{def:generic} is considered~\cite{goemans2009approximating,jegelka2011online}. 
\begin{defn}[Generic approximation~\cite{goemans2009approximating,jegelka2011online}]
The function $\ftg_t: 2^V \mapsto \RR$ is a generic approximation of $f_t$ defined as:
$
\ftg_t(S) = \sqrt{\sum_{i \in S} c_i},
$
for some $\mathbf{c} \in \RR^n$, and satisfies
$
\left(\ftg_t(S)\right)^2 \leq \left(f_t(S)\right)^2 \leq \gamma^2 \left(\ftg_t(S)\right)^2,
$
for all $S \subseteq 2^V$ and some $\gamma > 0$.
\label{def:generic}
\end{defn}
Interested readers are referred to~\cite{goemans2009approximating} for details about the constant~$\mathbf{c}$. 
We now consider the online submodular generic greedy algorithm (\texttt{OSGGA}), i.e., the generic approximation-based \texttt{OSGA}. \texttt{OSGGA} uses the following update:
\begin{equation}
S_t \in \argmin_{S \in \mathcal{S}} \left(\ftg_{t-1}(S)\right)^2,
\label{eq:update_generic}
\end{equation}
The update rule~\eqref{eq:update_generic} is equivalent to solving a mixed-integer program (MIP) with a linear objective function~\cite{jegelka2011online} and can, therefore, be solved efficiently using off-the-shelf solvers for feasible sets that are linear or convex if relaxed.

For the next result, we utilize the variation term $\tilde{V}_T^\text{g} = \sum_{t=2}^T \sqrt{\card\left( S_{t}^{\text{g},\star} \ominus S_{t-1}^{\text{g},\star}\right)}$, where $S_t^{\text{g},\star} \in \argmin_{S \in \mathcal{S}} (\ftg_{t}(S))^2$ can be effectively computed via MIP for linear or convex (if relaxed) problems. Similarly to \texttt{OSGA}, $\tilde{V}_T^\text{g}$ is based only on the optima of tractable problems. Let $\nu \geq \min_{t,S\in\mathcal{S}} f_t(S)> 0$, be lower bound on all round minima. We remark that the squared generic approximation function satisfies Assumption~\ref{ass:pseudo_L} with modulus $\tilde{L}^\text{g}$ because it is linear and bounded by Assumption~\ref{ass:bounded}. The $\alpha$-regret for the \texttt{OSGGA} is presented below. 

\begin{cor}
Suppose $\ftg_t$ is a generic approximation of $f_t$ which can be solved via MIP. Then the $\alpha$-regret of \texttt{OSGGA} is bounded above by
\begin{align*}
\textup{\R}^\textup{\text{d}}_\alpha(T) &\leq \frac{4 \alpha^2 \tilde{L}^\text{g}L}{(1+\alpha)\nu} \sum^T_{t=2} \sqrt{\card\left( S_{t}^{\text{g},\star} \ominus S_{t-1}^{\text{g},\star}\right)}\\
&= \frac{4 \alpha^2 \tilde{L}^\text{g}L}{(1+\alpha)\nu} \tilde{V}_T^\text{g},
\end{align*}
and is sublinear for $\tilde{V}_T^\text{g} < O(T)$.
\end{cor}
\begin{pf}
We based our proof on~\cite[Theorem 2 and Lemma 3]{jegelka2011online}. The $\alpha$-regret for update~\eqref{eq:update_generic} is
\begin{align*}
\R_\alpha^\text{d}(T) &= \sum_{t=1}^T f_t(S_t) - \alpha f_t(S_t^\star) \\
&=\sum_{t=1}^T \frac{\left(f_t(S_t) - \alpha f_t(S_t^\star)\right) \left(f_t(S_t) + \alpha f_t(S_t^\star)\right)}{\left(f_t(S_t) + \alpha f_t(S_t^\star)\right)}\\
&=\sum_{t=1}^T \frac{\left(f_t(S_t)\right)^2 - \alpha^2 \left(f_t(S_t^\star)\right)^2}{\left(f_t(S_t) + \alpha f_t(S_t^\star)\right)}
\end{align*}
\begin{align*}
\phantom{\R_\alpha^\text{d}(T)}
&\leq\sum_{t=1}^T \frac{\left(f_t(S_t)\right)^2 - \alpha^2 \left(f_t(S_t^\star)\right)^2}{(1+\alpha)\nu},
\end{align*}
where $\nu \geq \min_{t,S_t} f_t(S_t)> 0$. By Definition~\ref{def:generic}, we have
\begin{align*}
\R_\alpha^\text{d}(T) &\leq \frac{1}{(1+\alpha)\nu} \sum_{t=1}^T \gamma^2 \left(\tilde{f}_t^\text{g}(S_t)\right)^2 - \alpha^2 \left(\tilde{f}_t^\text{g}(S_t^\star)\right)^2\\
&\leq \frac{1}{(1+\alpha)\nu} \sum_{t=1}^T \gamma^2 \left(\tilde{f}_t^\text{g}(S_t)\right)^2 - \alpha^2 \left(\tilde{f}_t^\text{g}(\tilde{S}_t^{\text{g},\star})\right)^2,
\end{align*}
where $\tilde{S}_t^{\text{g},\star} \in \argmin_{S \in \SS} \tilde{f}_t^\text{g}(S)$. Using the update rule~\eqref{def:generic}, we obtain
\begin{align}
\R_\alpha^\text{d}(T) &\leq \frac{\gamma^2}{(1+\alpha)\nu} \sum_{t=1}^T  \left(\tilde{f}_t^\text{g}(\tilde{S}_{t-1}^{\text{g},\star})\right)^2 - \left(\tilde{f}_t^g(\tilde{S}_t^{\text{g},\star})\right)^2. \label{eq:inter_preL}
\end{align}
Thus, the $\alpha$-regret can be re-expressed as
\begin{align*}
\R_\alpha^\text{d}(T) &\leq \frac{\gamma^2}{(1+\alpha)\nu} \R^\text{d}\left(\left(\tilde{f}_t^\text{g}(S)\right)^2,T\right),
\end{align*}
where $\R^\text{d}\left(\left(\tilde{f}_t^g(S)\right)^2,T\right)$ is the ($1$-)regret of update~\eqref{eq:update_generic} when used on the problem $\min_{S \in \SS} \left(\tilde{f}_t^g(S)\right)^2$.

Using Theorem~\ref{eq:thm_greedy} with $\alpha=\beta=1$ in~\eqref{eq:inter_preL} yields
\begin{align*}
\R_\alpha^\text{d}(T) &\leq \frac{\gamma^2 \tilde{L}^\text{g}L}{\nu} \sum_{t=1}^T  \sqrt{\card\left(\tilde{S}_{t-1}^{\text{g},\star} \ominus \tilde{S}_t^{\text{g},\star} \right)} = \frac{\gamma^2 \tilde{L}^\text{g}L}{\nu} \tilde{V}^\text{g}_T,
\end{align*}
which completes the proof. $\hfill \qed$
\end{pf}
Hence, \texttt{OSGGA} leads to an $\alpha$-regret bound that has a form similar to Theorem~\ref{eq:thm_greedy}'s. 
Comparing \texttt{OSGGA} to \texttt{OSGA}, different variation terms and constant factors in the regret bounds are used, which can lead to an increase in its value. However, the former can always be used.

\section{Projected gradient descent-based approach}
\label{sec:projected}
In this section, we consider a convex optimization-based update to solve~\eqref{eq:prob_submodular}~\cite{jegelka2011online,hazan2012online}. Our approach leverages the Lov\'asz extension's convexity for submodular functions. We propose the online submodular projected gradient descent (\texttt{OSPGD}) based on the update defined as:
\begin{align}
\xx_{t+1} &= 
\proj_{\conv(\SS)} \xx_t - \eta \gg_t \label{eq:proj_1}\\
S_{t+1} &= \rounding_\SS(\xx_{t+1}), \label{eq:proj_2}
\end{align}
where $\eta>0$ is the descent stepsize, $\gg_t \in \partial \hf_t(\xx_t)$ is defined in~\eqref{eq:gradient}, $\proj_{\conv(\SS)}$ is the projection onto the convex hull of $\mathcal{S}$, and $\rounding_\SS(\xx_{t+1}):[0,1]^n \mapsto \SS$ is defined and exemplified in Section~\ref{ssec:sub}. \texttt{OSPGD} has the advantage over the greedy updates to be computationally very simple because it only performs a single projected gradient descent step. It requires only algebraic operations and a projection onto a convex set which can be readily computed in most instances, e.g., box constraints, making it amenable to large problems. The use of \texttt{OSPGD} is, however, limited because it requires a rounding algorithm which might not be available for all constrained problems.
For \texttt{OSPGD}, we extend the regret analysis of~\cite{jegelka2011online} and obtain the following regret bound.

\begin{thm} Suppose that a rounding algorithm with approximation guarantee $\alpha$ is used. Then, \texttt{OSPGD} with $\eta = \frac{\delta}{\sqrt{T}}$ leads to an $\alpha$-regret bounded from above by:
\begin{align*}
\textup{\R}^{\textup{d}}_\alpha(T) 
&= \alpha \left(\sqrt{n} \delta V_T + \frac{5n}{2\delta} + 4M\delta \right)\sqrt{T},
\end{align*}
and is sublinear if $V_T < O\left(\sqrt{T}\right)$.
\label{thm:LPGD_alpha}
\end{thm}
\begin{pf}
By definition, we have
\begin{align*}
\R^\text{d}_\alpha(T) &= \sum_{t=1}^T  f_t(S_t) - \alpha f_t(S_t^\star) \\
&\leq \sum_{t=1}^T  \alpha \hf_t(\xx_t) - \alpha f_t(S_t^\star),
\end{align*}
using the rounding algorithm approximation guarantee bound. Using the property of the Lova\'sz extension, we obtain
\begin{align}
\R^\text{d}_\alpha(T) \leq \alpha \sum_{t=1}^T  \hf_t(\xx_t) - \hf_t(\bm{\chi}_{S_t^\star}). \label{eq:inter_regret}
\end{align}
We then follow the standard proof techniques for the online gradient descent (\texttt{OGD}) from~\cite{zinkevich2003online,hazan2016introduction}. We have
\begin{align}
\left\| \xx_{t+1} - \cc_{S_t^\star} \right\|_2^2 &= \left\| \left( \proj_{\conv(\SS)} \xx_t - \eta \gg_t \right) - \cc_{S_t^\star} \right\|_2^2 \nonumber\\
&\leq\left\|  \xx_t - \eta \gg_t  - \cc_{S_t^\star} \right\|_2^2 \nonumber\\
&= \left\|  \xx_t - \cc_{S_t^\star} \right\|_2^2 - 2 \eta_t \gg_t^\top \left( \xx_t - \cc_{S_t^\star} \right)\nonumber\\
&\qquad + \eta^2 \left\| \gg_t \right\|_2^2 \nonumber\\
\Leftrightarrow \gg_t^\top \left( \xx_t - \cc_{S_t^\star} \right) &\leq \frac{1}{2 \eta} \left( \left\|  \xx_t - \cc_{S_t^\star} \right\|_2^2 \right.\label{eq:inter_LPGD}\\
&\qquad\left.- \left\| \xx_{t+1}- \cc_{S_t^\star} \right\|_2^2 \right)
+ \frac{\eta^2}{2} \left\| \gg_t \right\|_2^2 .\nonumber
\end{align}
The convexity of $\hf_t$ implies that for all $\xx, \mathbf{y} \in [0,1]^n$:
\begin{align*}
\hf_t(\xx) \geq \hf_t(\mathbf{y}) + \gg_t^\top (\xx - \mathbf{y}),
\end{align*}
for $\gg_t \in \partial \hf_t(\mathbf{y})$. Using $\xx = \cc_{S_t^\star}$ and $\mathbf{y} = \xx_t$, we obtain
\begin{equation}
\hf_t(\xx_t) -  \hf_t(\cc_{S_t^\star})  \leq \gg_t^\top (\cc_{S_t^\star} - \xx_t). \label{eq:inter_2_LPGD}
\end{equation}
Substituting~\eqref{eq:inter_LPGD}~and~\eqref{eq:inter_2_LPGD} in~\eqref{eq:inter_regret} leads to
\begin{align*}
\R^\text{d}_\alpha(T) &\leq \alpha \sum_{t=1}^T \frac{1}{2 \eta} \left( \left\|  \xx_t - \cc_{S_t^\star} \right\|_2^2 - \left\| \xx_{t+1} - \cc_{S_t^\star} \right\|_2^2 \right)\\ 
&\qquad + \alpha \sum_{t=1}^T \frac{\eta^2}{2} \left\| \gg_t \right\|_2^2\\
&= \alpha \sum_{t=1}^T \frac{1}{2\eta} \left( \left\| \xx_t \right\|_2^2 - \left\| \xx_{t+1} \right\|_2^2 \right) + \alpha \sum_{t=1}^T \frac{\eta}{2} \left\|\gg_t \right\|_2^2\\
& \hspace{-0.3cm}+ \alpha \sum_{t=2}^T \frac{1}{\eta} \xx_t^\top \left( \cc_{S_t^\star} - \cc_{S_{t+1}^\star}\right) - \xx_1^\top \cc_{S_1^\star} + \xx_{T+1}^\top \cc_{S_T^\star}
\end{align*}
\begin{align*}
&\leq \frac{\alpha}{2\eta} \left\| \xx_1 \right\|_2^2 
-\frac{\alpha}{\eta} \xx_1^\top \cc_{S_1^\star} + \frac{\alpha}{\eta} \xx_{T+1}^\top \cc_{S_{T}^\star}\\
&\quad + \frac{\alpha}{\eta} \sum_{t=2} \xx_t^\top \left( \cc_{S_t^\star} - \cc_{S_{t-1}^\star} \right) + \alpha \sum_{t=1}^T \frac{\eta}{2} \left\|\gg_t \right\|_2^2,
\end{align*}
where we have evaluated the telescoping sums to obtain the last line~\cite{hazan2016introduction}. Using~\cite[Lemma 1]{jegelka2011online}, we have $\left\|\gg_t\right\|_2 \leq 4 M$. The regret becomes
\[
\R^\text{d}_\alpha(T) \leq \frac{5n\alpha}{2\eta} + 4 \alpha M \eta T + \frac{\alpha \sqrt{n} }{\eta} \sum_{t=2}^T \left\| \cc_{S_t^\star} - \cc_{S_{t-1}^\star} \right\|_2,
\]
where we also used the fact that $\left\|\xx_t \right\|_2 \leq \sqrt{n}$. We remark that for a submodular function and its Lova\'sz extension pair, $\left\| \cc_{S_t^\star} - \cc_{S_{t-1}^\star} \right\|_2$ is equivalent to $\sqrt{\card (S_t^\star \ominus S_{t-1}^\star )}$. We now have
\[
\R^\text{d}_\alpha(T) \leq \frac{5n\alpha}{2\eta} + 4 \alpha M \eta T + \frac{\alpha \sqrt{n} }{\eta} \sum_{t=2}^T \sqrt{\card \left(S_t^\star \ominus S_{t-1}^\star\right)}.
\]
Setting $\eta = \frac{\delta}{\sqrt{T}}$, $\delta \in \RR_{>0}$ completes the proof. $\hfill \qed$
\end{pf}
In sum, we obtain an $\alpha$-regret bound that is of the same order w.r.t. $T$ and the variation of optima term as the standard online gradient descent for OCO problems~\cite{zinkevich2003online}, i.e., $O\left(\sqrt{T}\left(1+V_T\right)\right)$. Theorem~\ref{thm:LPGD_alpha}'s bound differs from~\cite{lesage2021online}'s as it holds asymptotically and admits constrained decision-making. This difference can be explained in part by the fact that the rounding algorithm in~\cite{lesage2021online} yields the stronger property $\hf_t(\xx_t) = \E\left[ f_t(S_t)\right]$ instead of $\xx_t = \E\left[ S_t\right]$ in our case because the Lova\'sz extension is linear in $\xx$.
Lastly, in comparison to Section~\ref{sec:greedy}'s approaches, we have traded higher algorithmic simplicity for a regret bound that now depends on both $T$ and a variation term, the latter of which needs to be less than $O(\sqrt{T})$ to ensure a sublinear bound. Except in the unconstrained case discussed later, the bound can be hard to characterize because it is a function of~\eqref{eq:prob_submodular}'s optima.

Lastly, if a randomized rounding technique is used to convert a continuous decision vector to a binary one, expected and high-probability regret bounds, i.e., where $\alpha=1$ as opposed to previous results, can be derived.
\begin{cor}
\label{cor:projected_random}
Consider a random rounding technique $\rounding_\SS: \conv(\SS) \mapsto \SS$ such that for $S = \rounding_\SS(\xx)$ we have $\E\left[ f_t(S)) \right] = \hat{f}_t(\xx)$. The expected and high-probability dynamic regret for \texttt{OSPGD} with $\frac{\delta}{\sqrt{T}}$ are bounded from above:
\begin{align*}
\E \left[\textup{\R}^\textup{\text{d}}(T)\right] &\leq \sqrt{nT}\delta V_T + \left( \frac{5n}{2\delta} + 4M\delta \right)\sqrt{T},\\
\textup{\R}^{\textup{d}}(T) &\leq \sqrt{nT} \delta V_T + \left( \frac{5n}{2\delta} + 4M\delta + 2M\delta\log\frac{1}{\epsilon}\right)\sqrt{T},
\end{align*}
with probability of at least $1-\epsilon$.
\end{cor}
\begin{pf}
We adapt Theorem~\ref{thm:LPGD_alpha}'s and~\cite[Theorem 1]{hazan2012online}'s proofs to the dynamic setting. First, for the expected bound, we have
\begin{align}
\E\left[ \R^\text{d}(T) \right] &= \sum_{t=1}^T  \E\left[f_t(S_t) \right] - \E\left[f_t(S_t^\star)\right] \label{eq:inter_RHS}\\
&= \sum_{t=1}^T  \hf_t(\xx_t) - f_t(S_t^\star) \nonumber\\
&= \sum_{t=1}^T  \hf_t(\xx_t) - \hf_t(\cc_{S_t^\star}). \nonumber
\end{align}
The bound then follows from Theorem~\ref{thm:LPGD_alpha}. Second, for the high probability bound, we use H\oe ffding inequality~\cite[Theorem 13]{hazan2012online}. With a probability of a least $1-\epsilon$, we have
\begin{equation}
\sum_{t=1} f_t(S_t) \leq \sum_{t=1}^T \E[f_t(S_t)] + M \sqrt{2T \log \frac{1}{\epsilon}}. \label{eq:inter_hoeffding}
\end{equation}
Substituting~\eqref{eq:inter_hoeffding} in the regret definition, we obtain
\begin{align}
\R^\text{d}{(T)} \leq \sum_{t=1}^T \left(\E[f_t(S_t)] - f_t(S_t^\star) \right) + M\sqrt{2 T \log \frac{1}{\epsilon}} .\label{eq:after_Hoeff}
\end{align}
We observe that the first term of~\eqref{eq:after_Hoeff}'s right-hand side and~\eqref{eq:inter_RHS}'s are identical. Using Theorem~\ref{thm:LPGD_alpha} with $\eta=\frac{\delta}{\sqrt{T}}$ in~\eqref{eq:after_Hoeff} yields the high probability regret bound. $\hfill \qed$
\end{pf}

For unconstrained problems, we have $\SS = 2^V$ and $\conv \SS = \left[0,1\right]^n$. Then, Corollary~\ref{cor:projected_random} holds when the randomized rounding procedure described in Section~\ref{ssec:sub} is implemented in \texttt{OSPGD}. At this time, $V_T$ can be calculated efficiently, if desired.

\section{Applications to electric power systems}
\label{sec:appli}
We apply \texttt{OSPGD} and \texttt{OSGA} to power system problems.

\subsection{Demand response for frequency regulation}
In demand response, a load aggregator is contracted by the system operator~\cite{siano2014demand}. The aggregator's mandate is to modulate the load power consumption to help out the grid, e.g., to mitigate renewable intermittency or reduce peak demand. Specifically, we consider frequency regulation services~\cite{mathieu2012state,taylor2016power,lesage2018setpoint}, i.e.,  load balancing on a fast timescale, e.g., 4 seconds. Advantages of demand response over other frequency regulation approaches, like battery energy storage and fast-ramping fuel-burning generation, include low deployment costs and sustainability~\cite{taylor2016power}.

\subsubsection{Setting}
Consider $N$ thermostatically controlled loads (TCLs), e.g., residential loads equipped with electric water heaters, heaters, or air conditioners, enrolled in the demand response program. Consider a program of duration $T$ in which decision rounds are indexed by $t$. Let $p_{n,t} \geq 0$ and $\Tilde{p}_{n,t} \geq 0$ be the power consumption of TCL $n \in \left\{1,2,\ldots, N\right\}$ when the load is flexible and inflexible, respectively. This formulation is similar to~\cite{lesage2021online}'s. Each load must stay in an acceptable temperature range, e.g., $\pm 0.5^\circ$C of the desired user temperature, to be flexible, i.e., to be controlled according to the aggregator's need. If the load temperature is too high or too low, the backup controller forces the load to be active or inactive accordingly, and its power consumption must be accounted for. 

At time $t$, the aggregator's objective is to track a regulation setpoint $r_t>0$ provided by the system operator by adjusting the TCL power consumption. 
In this work, we consider a setting in which the aggregator wants to deploy the minimum number of flexible loads such that the regulation signal is met. 
This problem can be formulated as an online dynamic submodular optimization problem using the objective function $f_t^\text{DR}: 2^V \mapsto \RR$,
\begin{equation}
\begin{aligned} 
f_t^\text{DR}(S_t)= &\sum_{A \subseteq 2^V} 
\left[  \left(\sum_{n \in A} u_{n,t} \right)^2 -  \left(\sum_{n \in V} u_{n,t} \right)^2\right] \cdot
\label{eq:obj_dr}\\
& \qquad \max\{0, |S_t \cap A| - |S_t \cup A| + 1 \} \mathbb{I}_{A \subseteq \mathcal{R}_t},
\end{aligned}\end{equation}
where $u_{n,t} = p_{n,t} + \Tilde{p}_{n,t}$, $\mathbb{I}$ is the indicator function which returns $1$ if the subscript is true and $0$ otherwise, and 
\[
\Rt_t = \left\{ S \subseteq 2^V \left| \sum_{n \in S} u_{n,t} \geq r_t \right.\right\}.
\]
In~\eqref{eq:obj_dr}, the term between brackets promotes partitions $A \subseteq 2^V$ with lower aggregated power. Then, the maximum term identifies to which partition $A$ the set $S$ belongs to, because it is equal to one if and only if $S = A$.  Lastly, the indicator function ensures that the set of dispatched loads is at least equal to the regulation signal. 

We apply \texttt{OSPGD} to this problem. In terms of standard online optimization, this corresponds to a quadratic program with time-dependent binary constraints, which, to this day, has not been investigated. To the authors' best knowledge, no other approach has been shown to have provable performance in this context.

We randomly generate loads' parameters similarly to~\cite{mathieu2015tcl}. We use the same constraints and logical rules as~\cite{lesage2021online} and omit the lockout constraint. The thermodynamic model is based on~\cite{mathieu2012state}. We consider TCLs equipped with air conditioners. 

\subsubsection{Numerical results}
We deploy 15 TCLs to track a vanishing sinusoidal regulation signal subject to Perlin noise~\cite{perlin2002}. We compare \texttt{OSPGD} to the closest work to ours, \texttt{bOGD}~\cite{lesage2021online}. We note that, in this setting, \texttt{bOGD}'s regret analysis does not hold. Lastly, we provide the round optimum, which we denote $\texttt{OSPGD}^\star_t$.

\begin{figure}[tb]
  \begin{subfigure}[t]{.48\textwidth}
    \centering
    \includegraphics[trim=0 0 0.2cm 0,clip,width=\linewidth]{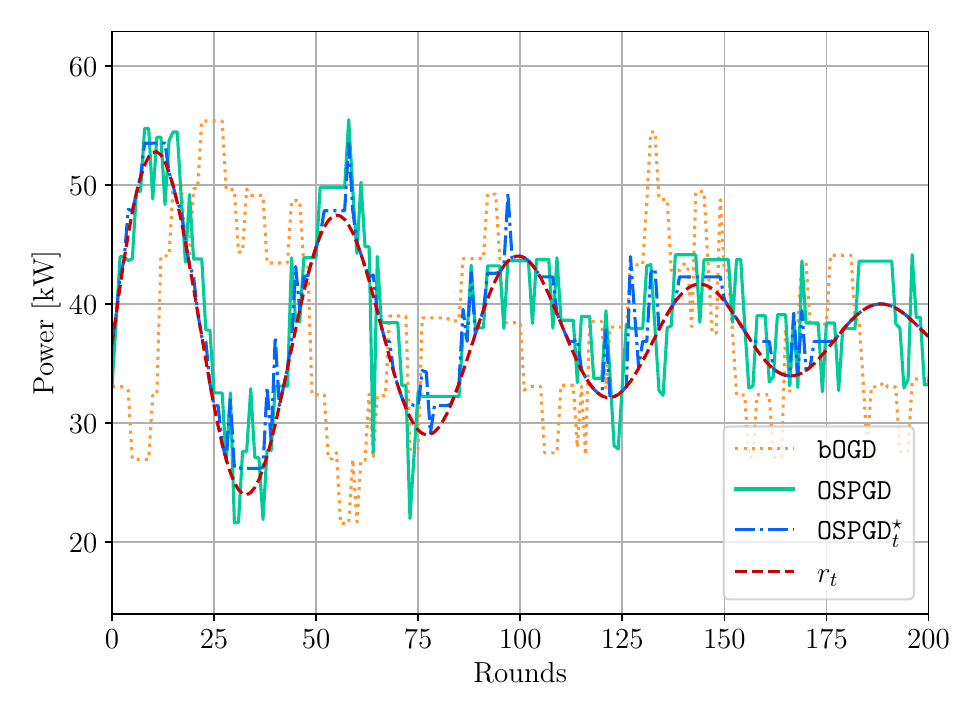}
    \caption{Setpoint tracking}
    \label{subfig:DRtracking}
  \end{subfigure}
  \hfill
  \begin{subfigure}[tb]{.48\textwidth}
    \centering
    \includegraphics[trim=0 0 0.2cm 0,clip,width=\linewidth]{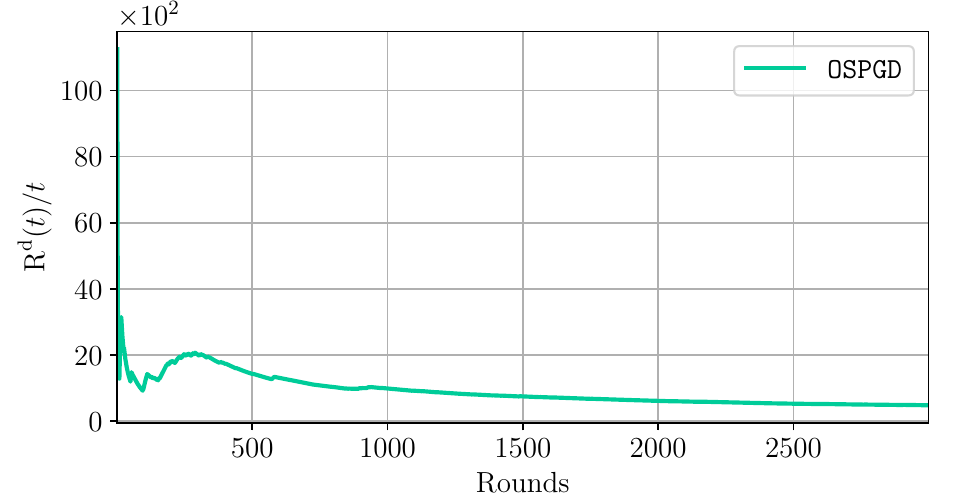}
    \caption{Time-averaged experimental dynamic regret for \texttt{OSPGD}}
     \label{subfig:DRregret}
  \end{subfigure}
  \caption{Demand response with 15 loads}
  \label{fig:DR}
\end{figure}

As shown in~Figure~\ref{subfig:DRtracking},  \texttt{OSPGD} outperforms \texttt{bOGD} and offers good setpoint tracking. The tracking root-mean-square error (RMSE) over $3000$ rounds in this case is $1.406$~kW for $\texttt{OSPGD}^\star_t$, $3.227$~kW for \texttt{OSPGD}, and $6.928$~kW for \texttt{bOGD}. Figure~\ref{subfig:DRregret} presents \texttt{OSPGD}'s time-averaged dynamic $1$-regret. The vanishing time-averaged regret implies a sublinear regret.

\subsection{Real-Time Network Reconfiguration}
Electric distribution networks generally possess a radial topology~\cite{taylor2015convex}. Their topology is controlled via switches located throughout the network. By opening and closing different switches, the topology can be modified, for example, to minimize active power losses or line congestion, and, thus, to increase the grid efficiency~\cite{rao2012power,mishra2017comprehensive,khodabakhsh2017submodular}. The set of switch statuses must always induce a radial network topology while assuring that all loads are supplied.

Distribution grids with high penetration of grid-edge/behind-the-meter technologies~\cite{paudyal2020impact}, e.g., electric vehicles, residential solar panels, or demand response, can experience large, fast-ramping variations in power demand at the different buses. These rapid changes in loading are out of the distribution system operator's control and can lead to network perturbations, e.g., over/under-voltage, line congestion, etc.~\cite{leou2013stochastic,navarro2015probabilistic}. To mitigate incidents, the system operator can preemptively configure the distribution network by altering its topology. Remotely-activated switches allow fast network reconfiguration (NR) and can be used to adapt to the load demand in real-time, viz., to prevent line congestion or to reduce active power losses.

\subsubsection{Setting}

We consider a distribution network consisting of a set of static powerlines $\overline{\mathcal{L}}$, a set of loads $\mathcal{N}$, and a set of lines equipped with switches $V=\left\{1,2,\ldots, \overline{S} \right\}, \ \overline{S}\in \mathbb{N}$. Let $\mathcal{L}(S_t) = \overline{\mathcal{L}} \cup S_t$ where $S_t \in 2^V$ be the set of all powerlines active at time $t$, i.e., the static line set augmented by the lines with closed switches $S_t$. The set $\mathcal{L}(S_t)$ is subject to two constraints; it must be such that (i) the network topology is radial and (ii) all loads are connected.

Let $p_{i,t} \geq 0$ and $q_{i,t} \geq 0$, be the active and reactive power demand, respectively, at bus $i \in \mathcal{N}$ and time $t$. Let $\mathcal{N}_r$ be the set of feeder nodes. Let $P_{ij,t} \in \RR$ and $Q_{ij,t} \in \RR$ be, respectively, the active and reactive power flowing from node $i$ to $j$ if $ij \in \mathcal{L}(S_t)$. Let $P_{ij,t} = Q_{ij,t} = 0$ if $ij \notin \mathcal{L}(S_t)$. Line $ij$'s apparent power is denoted by $A_{ij,t} = P_{ij,t} + j Q_{ij,t}$. Let $v_i \in \mathbb{C}$ be the voltage at node $i$,  $I_{ij}\in \mathbb{C}$ be the current flowing in line $ij$, and $y_{ij}\in \mathbb{C}$ be the admittance of line $ij$. Let notation $\overline{x}$ and $\underline{x}$ represent upper and lower bounds on any given parameter $x$.

To minimize active power losses in distribution grids, the NR problem can be cast using the objective function $f_t^\text{NR}: 2^V \mapsto \RR$ presented in~\eqref{eq:reconf} where power losses on line $ij$ at time $t$ are defined as $|v_{i,t} - v_{j,t}|^2 y_{ij}^*$. In~\eqref{eq:reconf}, the spanning tree constraint ensures that the network topology is radial and connects all loads to the source node. The other constraint ensure that the power flow (PF), which models the electric network's physics, respects all operational constraints while meeting power demand.
\begin{align}
& \min_{S_t \subseteq 2^V} & & f_t^{\text{NR}}(S) = {\sum_{ij \in \mathcal{L}(S)} |v_{i,t} - v_{j,t}|^2 y_{ij}^*}\nonumber\\
& \text{subject to}
& & \mathcal{L}(S_t) \subseteq \textrm{SpanningTree($\mathcal{N}$)}\label{eq:reconf}\\
& & &
\hspace{-1.25cm}\begin{Bmatrix}
\{v_{i,t}\}_{i\in\mathcal{N}}\\
\{P_{ij,t}, Q_{ij,t}\}_{ij \in \mathcal{L}(S_t)}
\end{Bmatrix}
\in \textrm{PF}(\{p_{i,t}, q_{i,t}\}_{i \in \mathcal{N}}, \mathcal{L}(S_t)). \nonumber
\end{align}

\subsubsection{Weakly-meshed Approximation}
Finding the optimal configuration of a radial network is NP-hard. We re-express (\ref{eq:reconf}) as an online dynamic submodular optimization problem, which can then be solved in real-time.

When the radiality constraint is relaxed, the network, in which the set of active powerlines is $\overline{\mathcal{L}} \cup V$, referred to as the weakly-meshed network (WMN), is a good solution, if not optimal, for loss minimization~\cite{AHMADI2015MST}. Using the WMN as a starting point, our goal is to find the radial network that best imitates its power flow. This can be done by first computing the WMN power flow. Then, a minimum spanning tree (MST) algorithm (e.g., Prim’s algorithm~\cite{prim1957}) with edge weights set as the negative line currents $-I_{ij}$ obtained from the power flow, is used. The MST is fast and guarantees radiality. By removing the edges with lower currents, the MST returns a radial network with a power flow pattern similar to the WMN as demonstrated by~\cite{AHMADI2015MST}. We note that in all evaluations, the resulting topology admitted a feasible power flow with respect to the original AC power flow constraints. If infeasible, the resulting topology could be projected onto the set induced by these constraints. Finally, we can approximate (\ref{eq:reconf}) by the following online dynamic submodular problem: 
\begin{equation}
\begin{aligned}
& \min_{{S_t} \in 2^V} & & f_t^\text{WM}({S_t}) = \sum_{ij \in {\overline{\mathcal{L}}}} I_{ij,t} + \sum_{ij \in {\mathcal{L}(S_t)}} -I_{ij,t}\\
& \text{subject to}
& & \mathcal{L}({S_t}) \subseteq \text{SpanningTree($\mathcal{N}$)},
\end{aligned}\label{eq:wmreconf}
\end{equation}
where $I_{ij,t}$ is an online parameter extracted from power flow computations, e.g.~\cite{pandapower.2018}, of the WMN.

Because (\ref{eq:wmreconf}) is submodular, and can be solved to optimality in polynomial time using a MST algorithm like Prim's~\cite{prim1957} over the WMN, we apply our \texttt{OSGA} for online reconfiguration. The process is summarized in Algorithm \ref{alg:NR}.
In Algorithm \ref{alg:NR}, $M$ is a large constant. We remark that in the case of multiple feeders, we temporarily add virtual lines between the different sources (generators) in the MST algorithm to ensure radiality, see steps 5$-$6. These lines are then removed from $S_t$, see step 8.

\begin{algorithm}[tb]
\caption{\texttt{OSGA} for Real-time NR}
\begin{algorithmic}[1]
\FOR {$t = 1,2,...,T$} 
\STATE Reconfigure the network according to $S_t$.
\STATE Observe new online parameters: ${p}_{i,t}$, $q_{i,t}, \forall i  \in  \mathcal{N}$.
\STATE Compute the power flow of the WMN with the Newton-Raphson algorithm and extract $I_{ij,t} \ \forall ij \in \overline{\mathcal{L}} \cup V$.
\IF {$\text{card}(\mathcal{N}_r) > 1$}
\STATE Set $I_{ij,t} = M, \ \forall ij \text{ where } \{i,j \in \mathcal{N}_r \text{ and } i \neq j$\}
\ENDIF
\STATE Update $S_{t+1}$ via a greedy MST algorithm: 
\begin{equation*}
\begin{aligned}
\Tilde{S}_{t+1} &= \argmin \: f_t^{\text{WM}}({S})\\
&\quad\text{subject to} \:\: \mathcal{L}({S}) \subseteq \text{SpanningTree($\mathcal{N}$)}\\
S_{t+1} &\in \Tilde{S}_{t+1} \setminus \{\text{ virtual lines }\}
\end{aligned}
\end{equation*}
\ENDFOR
\end{algorithmic}
\label{alg:NR}
\end{algorithm}

\subsubsection{Numerical Results}
\begin{figure}[tb]
    \centering
    \includegraphics[trim=0 0 0.5cm 0,clip, width=0.48\textwidth, height=0.5\linewidth]{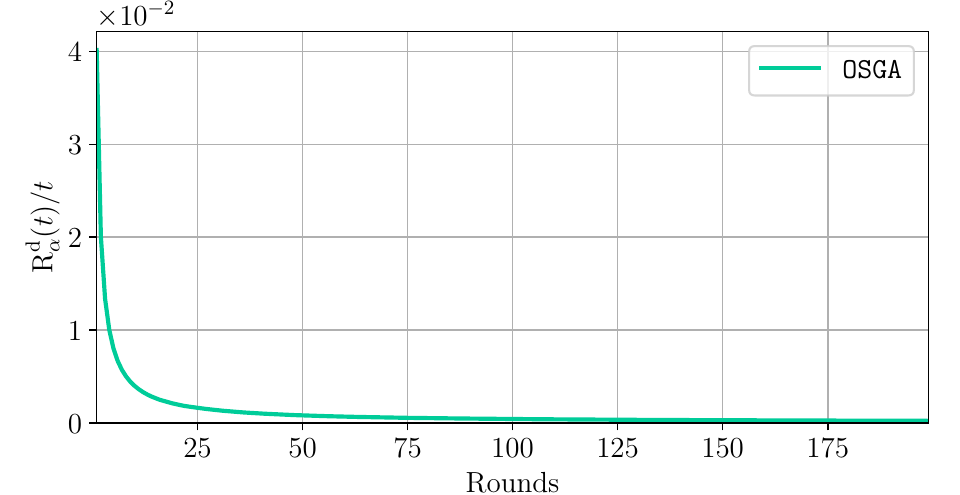}
    \caption{Time-averaged experimental dynamic regret for \texttt{OSGA}}
    \label{fig:tregret}
\end{figure}

For this section, we consider the IEEE standardized 33-bus/1-feeder (33b/1f)~\cite{Baran1989_33bus} and the 135-bus/8-feeder (135b/8f)~\cite{reds} distribution networks with the added modification, on both networks, that every line is equipped with a switch to fully benefit from the flexibility of online optimization. At each round, we add randomly generated Perlin noise~\cite{perlin2002} on ${p}_{i,t}$, $q_{i,t}, \forall i \in \mathcal{N}$ to model uncertainty. Figure \ref{fig:tregret} illustrates \texttt{OSGA}'s sublinear dynamic $1$-regret. This is depicted by the vanishing time-averaged regret.

We now compare \texttt{OSGA} to its offline counterparts solved in hindsight both dynamically ($\texttt{OSGA}^\star_t$) and statically ($\texttt{OSGA}^\star$) over the time horizon. We note that hindsight solutions only serve analysis purposes and have no practical application. We benchmark our approach, in the simpler network (33b/1f), to the closest work in OCO (\texttt{bOGD})~\cite{lesage2021online} to which we must add a projection on the feasible power flow set to handle operational constraints of the grid. We also compare \texttt{OSGA} to a round-optimal offline configuration, with a limited flexibility of 9 switches, based on the second-order cone relaxation power flow (\texttt{SOCR-9SW})~\cite{taylor2015convex}, which require much more computational power. Lastly, we present the case where a random feasible reconfiguration is implemented each round.

\begin{figure}[tb]
    \centering
    \includegraphics[trim=0 0 0.5cm 0,clip,width=0.48\textwidth]{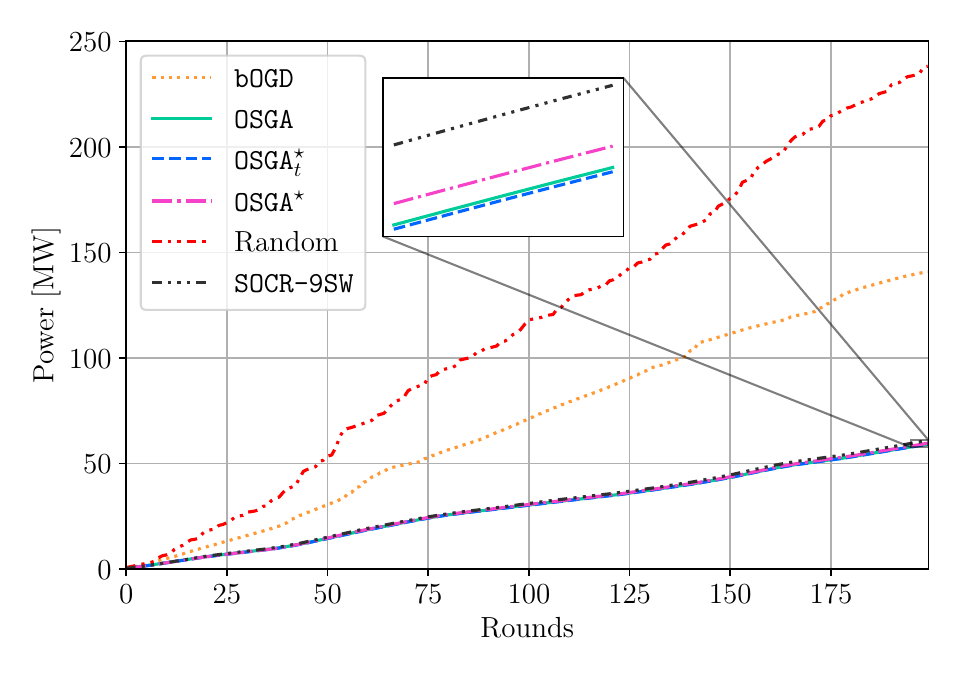}
    \caption{Cumulative power losses (33b/1f)}
    \label{fig:enter-label}
\end{figure}

\begin{figure}[tb]
   \begin{subfigure}[t]{.25\textwidth}
     \centering
     \includegraphics[trim=0.6cm 0 5cm 0.5cm,clip,width=\linewidth]{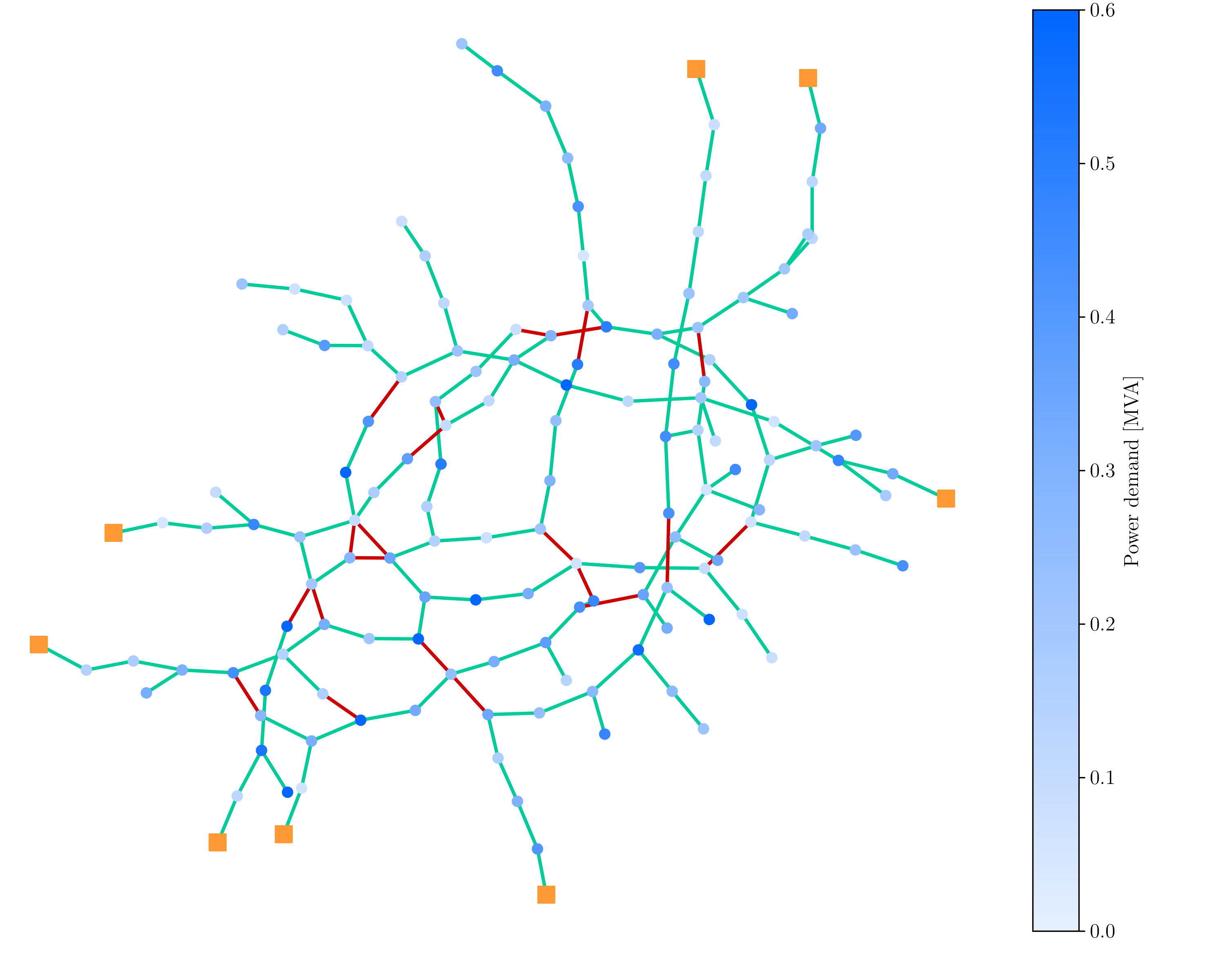}
     \caption{Round 20}
   \end{subfigure}%
   \begin{subfigure}[t]{.25\textwidth}
     \centering
     \includegraphics[trim=0.6cm 0 5cm 0.5cm,clip,width=\linewidth]{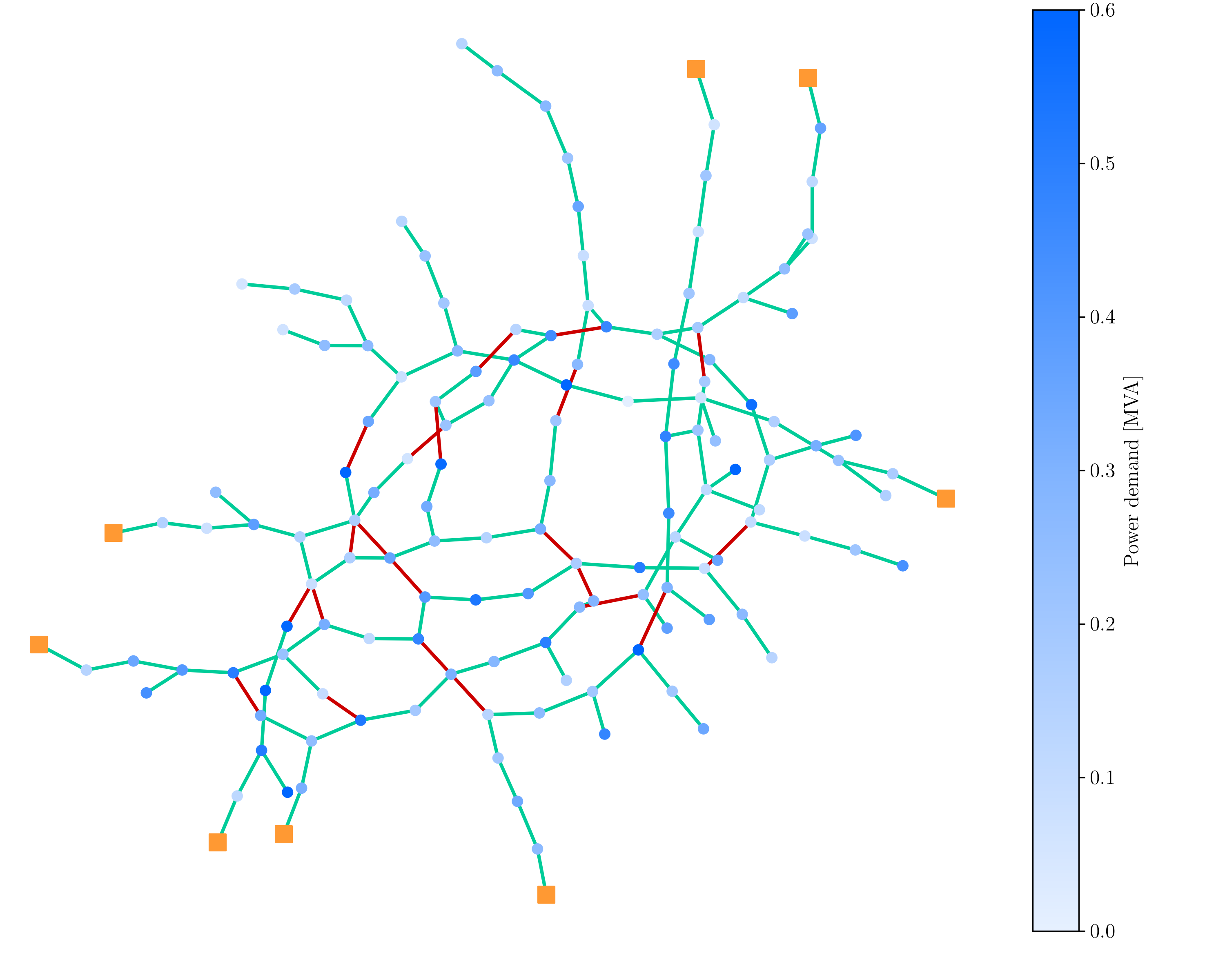}
     \caption{Round 21}
   \end{subfigure}
   \begin{subfigure}[t]{.25\textwidth}
     \centering
     \includegraphics[trim=0.6cm 0 5cm 0.5cm,clip,width=\linewidth]{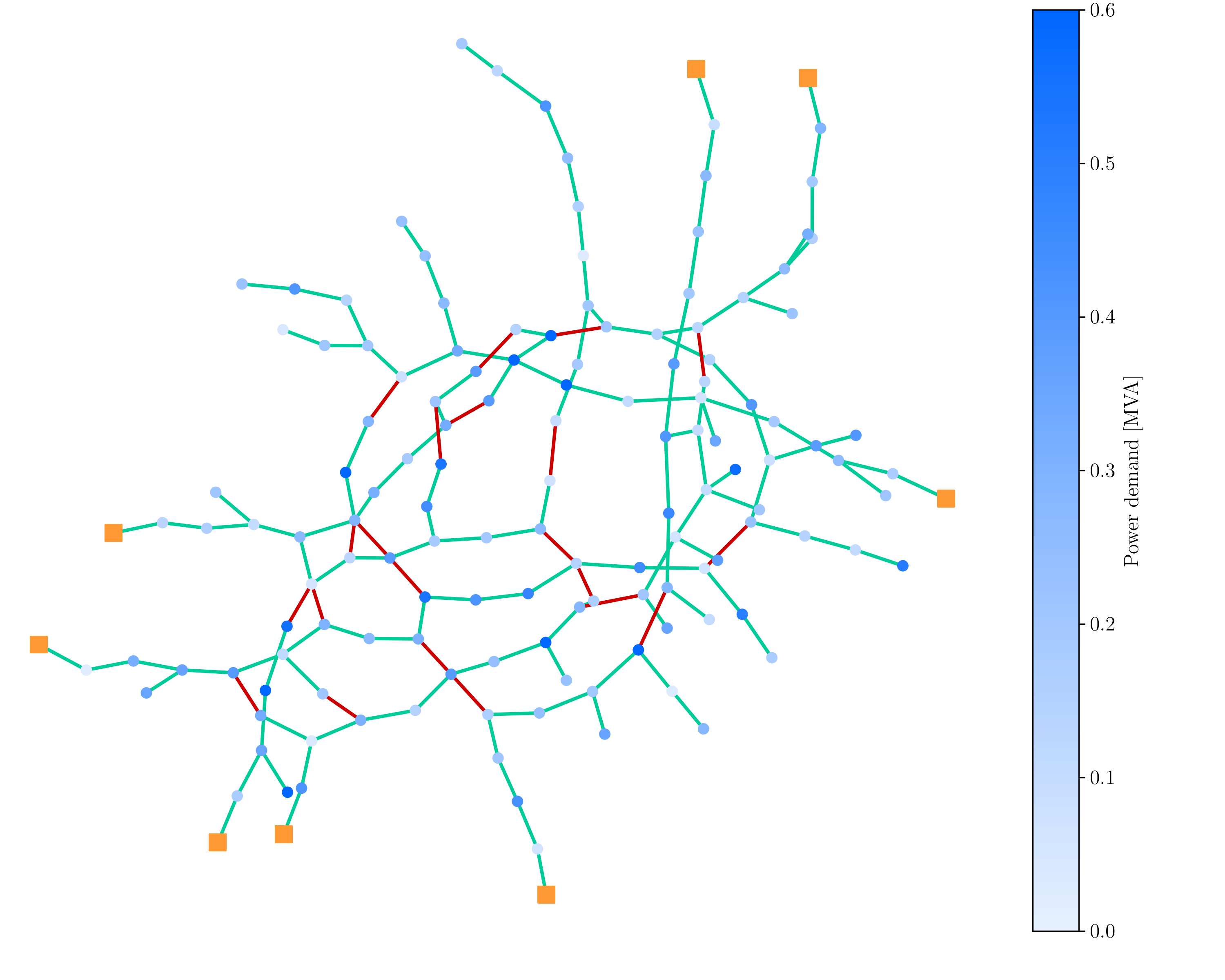}
     \caption{Round 22}
   \end{subfigure}%
   \begin{subfigure}[t]{.25\textwidth}
     \centering
     \includegraphics[trim=0.6cm 0 5cm 0.5cm,clip,width=\linewidth]{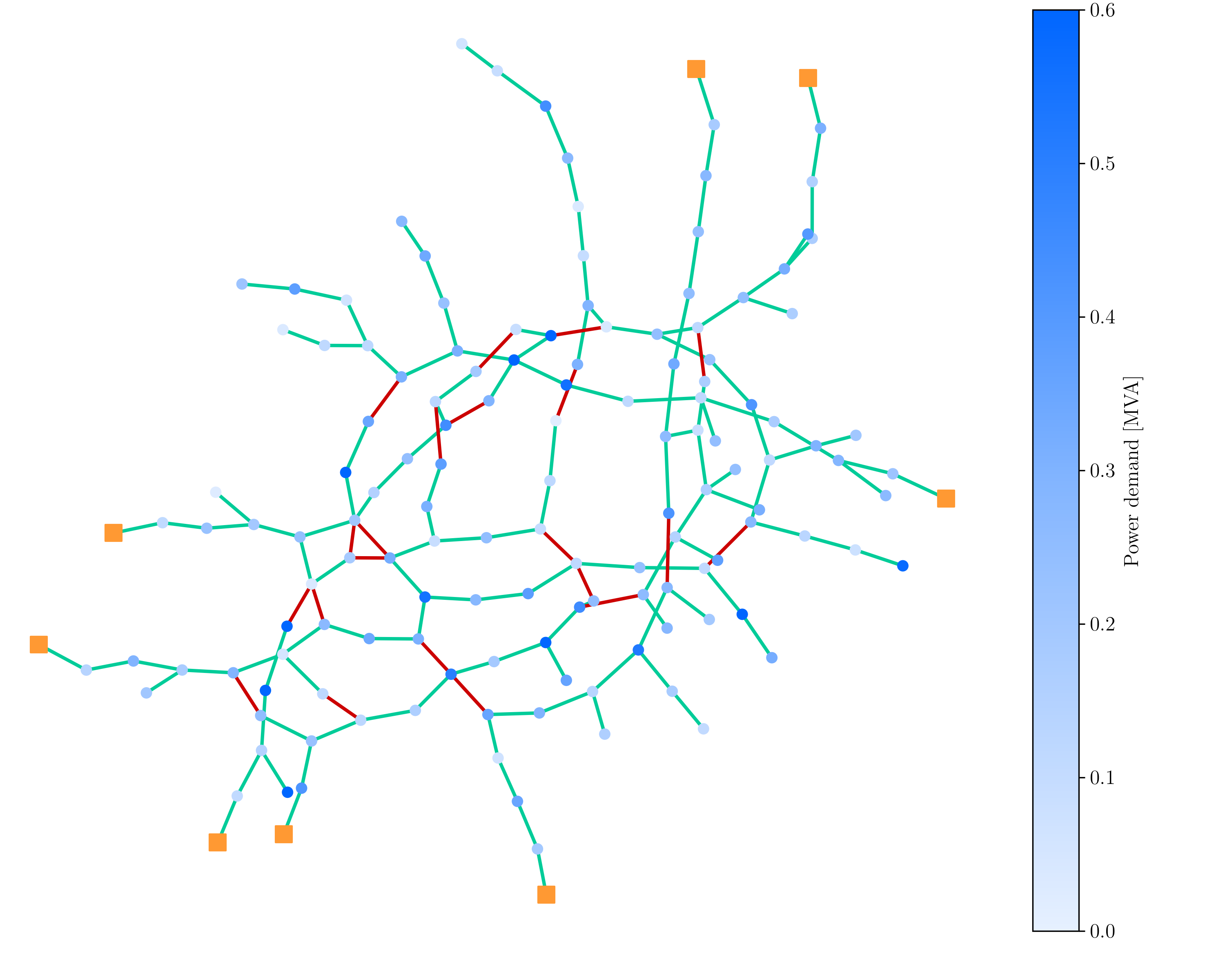}
     \caption{Round 23}
   \end{subfigure}
   \caption{Network reconfiguration (135b/8f)}
   \label{roundconfig}
 \end{figure}

Figure \ref{roundconfig} presents snapshots of the 135/8f NR at different rounds according to the apparent power demand at each node. The demand is represented by a light-dark scale: the darker the node the higher the demand is. Closed and open switches are pictured in green and red, respectively. Squares are generators. Radiality is always preserved.

In sum, \texttt{OSGA} performs systematically better than \texttt{bOGD} and is considerably faster because it stands on a fast MST heuristic. It also scales easily to bigger networks and guarantees radiality without the need for a projection step. \texttt{OSGA} also outperforms ${\texttt{OSGA}^\star}$ while maintaining a small performance gap with ${\texttt{OSGA}^\star_t}$, the round optimal solution computed in hindsight. For example, we observed a total power loss increase of $0.038\%$ for \texttt{OSGA} and of $1.545\%$ for ${\texttt{OSGA}^\star}$, after $400$ rounds, when compared to ${\texttt{OSGA}^\star_t}$ on 135b/8f in the simulation leading to~Figure~\ref{roundconfig}.

\section{Conclusion}
\label{sec:conclu}
In this work, we investigate online binary optimization in dynamic settings. We consider submodular objective functions and general binary constraints. 
We first assume an approximation of the objective function which can be minimized in polynomial time exists. We propose \texttt{OSGA} that solves the previous round approximation as a proxy and in doing so, circumvents the NP-hardness of the original problem. We adapt our approach to a generic but weaker approximation that can be used to recast general submodular problems in a simpler form. Second, aiming at algorithmic simplicity, we formulate \texttt{OSPGD} which leverages the Lov\'asz extension and convex optimization. For all our algorithms, we provide a dynamic regret analysis. We show that \texttt{OSGA} and \texttt{OSPGD} possess, respectively, a dynamic regret bound that is similar to the tightest bound w.r.t. the time horizon and the (tractable) round optimum variation in the literature and to the \texttt{OGD} used in online convex optimization. 

Finally, we present two applications of our approaches in electric power systems. First, \texttt{OSPGD} is employed to dispatch demand response resources, viz., thermostatic loads, to mitigate fast-timescale power imbalances. 
Second, \texttt{OSGA} is used to minimize active power losses in distribution networks via real-time reconfiguration, i.e., closing and opening switches in the network to better shape its topology.
 
Next, time-varying binary constraints, i.e., constraints that similarly to the objective function are observed only at the end of a round while needing to be satisfied in the long-run, and bandit feedback will be investigated. This will be done by combining, e.g., the idea behind~\texttt{OSPGD}, and a specialized approach like \texttt{MOSP}~\cite{chen2017online} and the point-wise gradient estimator from~\cite{hazan2016introduction}, respectively.

\begin{ack}                              
This work was funded by the Institute for Data Valorization (IVADO) and by the National Science and Engineering Research Council of Canada (NSERC).
\end{ack}

\bibliographystyle{plain}

\end{document}